\documentclass[12pt]{amsart}
\usepackage{mathrsfs}
\usepackage{graphicx}
\usepackage{enumerate}
\usepackage{amssymb,amsmath,amsthm, amscd}
\usepackage[all]{xy}
\usepackage{tikz}

\theoremstyle{plain}
\newtheorem{thm}{Theorem}[section]
\newtheorem{lem}[thm]{Lemma}

\newtheorem{cor}[thm]{Corollary}
\newtheorem{defn-lem}[thm]{Definition-Lemma}

\newtheorem{prop}[thm]{Proposition}
\newtheorem{question}[thm]{Question}

\theoremstyle{definition}
\newtheorem{defn}[thm]{Definition}
\newtheorem{ex}[thm]{Example}
\newtheorem{rem}[thm]{Remark}

\newcommand{\Z}{\mathbb{Z}}

\newcommand{\C}{\mathbb{C}}
\newcommand{\N}{\mathbb{N}}

\def\Index{\operatorname{Index}}

\def\id{\operatorname{id}}
\def\dist{\operatorname{dist}}
\def \Image{\operatorname{Image}}

\def\Ad{\operatorname{Ad}}
\def\Aut{\operatorname{Aut}}
\newcommand{\mc}{\mathcal}
\newcommand{\mb}{\mathbb}
\numberwithin{equation}{section}

\begin{document}
\title [Examples of inclusions with the Rokhlin property ]
       {The Rokhlin Inclusions with integer and Non-integer Index}

\begin{abstract}
We construct new examples of inclusions of unital $C\sp*$-algebras of index-finite type with the Rokhlin property motivated by a broader attempt to understand the range of such inclusions beyond known models.  In the course of this development, we observe an interesting phenomenon: inclusions with integer index, though not assumed to arise from group actions, exhibit internal behavior consistent with classical symmetry. In contrast, we construct inclusions with irrational index whose Rokhlin or tracial Rokhlin property arises from quantum symmetries - such as subfactor theory or more advanced tensor category action on Kirchberg algebras - and which cannot be modeled as fixed point algebras under any finite group action or  finite dimensional Hopf  $C\sp*$-algebra action. To our knowledge, these provide the first examples of inclusion with the  tracial Rokhlin property not arising from a finite group action.         
\end{abstract}

\author { Hyun Ho \, Lee}

\address {Department of Mathematics\\
         University of Ulsan\\
         Ulsan, 44610 South Korea \  \&
         School of Mathematics\\
         Korea institute of Advance Study\\
         Seoul, 02455 South Korea   }
\email{hadamard@ulsan.ac.kr}

\author{Hiroyuki Osaka}
\address{Department of Mathematical Sciences\\
Ritsumeikan Univeristy \\
Kusatsu, Shiaga, Japan }
\email{osaka@se.ritsumei.ac.jp}

\author{Tamotsu Teruya}
\address{Faculty of Education \\
Gunma University \\
Maebashi City, Gunma, Japan
}
\email{teruya@gunma-u.ac.jp}


\date{}
%\thanks{ }
\maketitle
\section{Introduction} 
The notion of the Rokhlin property first appeared as a fundamental concept in ergodic theory and was systematically introduced into $C^*$-algebra theory in the 1980s by R.H. Herman and A. Ocneanu for automorphisms of UHF $C\sp*$-algebras \cite{HO}. The property was subsequently generalized by A. Kishimoto \cite{K, K1, K2}, M. Izumi \cite{Izumi1, Izumi2}, and others \cite{AK, BSV, G, HJ}, extending its scope to a wider class of $C^*$-algebras and group actions. The Rokhlin property is a strong approximation condition that implies a given action behaves like a trivial action on a "large" part of the algebra. This property has proven to be a powerful tool, particularly in the classification of $C^*$-algebras, where it played a crucial role in the study of crossed product $C\sp*$-algebras.

\bigskip
This line of inquiry is situated within the broader context of the Elliott classification program, a monumental effort to classify a wide range of simple, separable, nuclear $C^*$-algebras using algebraic invariants such as K-theory. The success of this program hinges on the use of key technical notions that facilitate classification, including the concepts of strongly self-absorbing $C^*$-algebras (like the Jiang-Su algebra  $\mathcal{Z}$) \cite{TW}, and the use of approximation techniques such as approximate unitary equivalence and approximate intertwining \cite{Ell}. These tools provide the necessary machinery to show when two seemingly different $C^*$-algebras are, in fact, isomorphic. The Rokhlin property, or its variants, is a strong technical condition that has been instrumental in verifying these approximation conditions, proving that important regularity properties (e.g., $\mathcal{Z}$-stability, finite nuclear dimension) are preserved under group actions and inclusions \cite{BS, LO, OT1}.

\bigskip
While substantial progress has been achieved in fully understanding the Rokhlin property within the realm of group actions, leading to a rich and mature theory, the landscape remains notably less explored when we shift our focus to inclusions of $C\sp*$-algebras-particularly those inclusions that are not inherently generated by group actions. This represents a considerable gap in current knowledge. Indeed, a central yet implicit question has remained open: is there any inclusion of index-finite type which does not arise from classical group actions?  

\bigskip
Our investigation began with the goal of constructing Rokhlin inclusions with integer index that were genuinely independent of group actions. However, the initial attempts revealed the surprising robustness of  the classical case: the examples seemingly irrelevant of group actions, detailed in Section \ref{S:2}, invariably admitted an internal group-theoretic interpretation.  This pivotal finding led us to confront a deeper question: if the integer-index world is tightly. bound to classical symmetry, could genuinely non-classical Rokhlin phenomena be found elsewhere?  This paper provides an affirmative answer by constructing a new, extensive class of Rokhlin inclusions derived from quantum symmetries from Jones' subfactor theory to the more sophisticated action of unitary tensor categories. Our constructions are systematic in the sense that a quantum symmetry gives rise to  a conditional expectation on strongly self-absorbing $C\sp*$-algebras $A$, then we extend it to the infinite tensor product of $A$ leveraging key tools from the Elliott program,  the approximate unitary equivalence and the approximate intertwining via the inductive limit. This work thereby demonstrates the power of this toolkit beyond its original goal to show that two algebras are isomorphic, repurposing it for the explicit construction of new $C\sp*$-algebraic phenomena.        

\bigskip 
A key feature of our constructions is that, in stark contrast to the classical case, they naturally produce inclusions with non-integer indices, providing immediate and incontrovertible proof of their quantum origin. The power of our framework extends further, as the same quantum mechanisms are shown to generate inclusions with integer indices $\ge 4$. This reveals a more profound principle: the fundamental distinction is not merely a partition based on the integrality of the index, but rather between symmetries of classical origin(finite groups) and those of quantum origin(tensor categories). 

\bigskip
The primary contribution of this work is the application of this framework to solve a long-standing problem in the field. We construct the first known example of an inclusion of simple $C\sp*$-algebras possessing the tracial Rokhlin property that is provably not derived from a finite group action.  This result fundamentally expands our understanding of this key regularity property, demonstrating that it is not intrinsically limited to the classical context of finite group actions. By establishing a new and fertile source of Rokhlin-type phenomena, this work opens new avenues for the study of quantum symmetries within the broader $C\sp*$-algebra classification program.
\newpage 
\section{Rokhlin inclusions with integer index: Classical Patterns}\label{S:2}
The following examples are technically elementary and based on well-understood $C\sp*$-algebras.  
However, they have not been previously studied in the context of the Rokhlin property or constructed independently of group actions.  
Moreover, this paper represents the first systematic attempt to explore Rokhlin inclusions through direct construction, rather than those arising from finite group actions.  Although the examples here are conceptually accessible, their emergent internal symmetry—despite the absence of any assumed group action—offers a striking point of comparison with the more quantum examples in Section 3.  Thus, this section serves both as a foundational motivation and a structural baseline for what follows.

\bigskip
Throughout the paper, we only consider inclusions of unital $C\sp*$-algebras with finite index which was introduced by Y. Watatani \cite{Wa}. 
 
\begin{defn}[Y. Watatani]
Let $P\subset A$ be an inclusion of unital $C\sp*$-algebras with a conditional expectation $E: A \to P$,  which is a completely positive contraction  such that $E(b)=b$,  $E(xb)=E(x)b$ and $E(cx)=cE(x)$ for all $x\in A$ and $b, c \in P$ or shortly projection map of norm one \cite [II.6.10.1]{Bl}.  A quasi-basis for $E$ is a finite set of ordered pairs $\{(u_i, v_i) \}_{i=1}^n$ such that  for every $a \in A$ 
\[ a =\sum_{i=1}^n u_i E(v_i a) =\sum_{i=1}^n E(a u_i) v_i. \]
When $\{(u_i, v_i) \}_{i=1}^n$ is a quasi-basis for $E$, $\Index E$ is defined  by $\displaystyle \Index E =\sum_{i=1}^n u_iv_i$ which is independent of the choice of the quasi-basis.  In this case, it is said that the inclusion $P \subset A$ is of index-finite type \cite[Definition 1.2.2]{Wa}. 
\end{defn}

For a $C\sp*$-algebra $A$, recall that the sequence algebra $A^{\infty}= l^{\infty}(\mathbb{N}, A)/ c_0(A)$ where $c_0(A)$ is the ideal of sequences whose norm limit is $0$. We embed $A$ into $A^{\infty}$ by the equivalence classes of constant sequences and call $A^{\infty}\cap A'$ the central sequence algebra of $A$.  

\begin{defn}\cite[Definition 3.1]{OKT} \label{D:Rokhlin}
Let $P\subset A$ be an inclusion of unital $C\sp*$-algebras  of index-finite type. Then we say that $E$ has the Rokhlin property if there exist a nonzero projection $e$ in $A^{\infty}\cap A'$, the central sequence algebra of $A$,  satisfying the following two conditions; 
\begin{equation}
(\Index E) E^{\infty}(e)=1_{A^{\infty}}, 
\end{equation}  
and 
\begin{equation}
\text{the map}\, \, A \ni x \mapsto xe \in A^{\infty} \,\, \text{is injective.} 
\end{equation}
where $E^{\infty}$ is the conditional expectation from $A^{\infty}$ onto $P^{\infty}$ induced by $E$.  Moreover, when $E$ has the Rokhlin property, then we call $e$ the Rokhlin projection. 
\end{defn}

To investigate how the integer index value of an inclusion with the Rokhlin property indicates an intrinsic group structure, we begin with a commutative but non-simple case first.      

\begin{ex} \label{Ex:1}

Let $X_1, X_2, \dots, X_{n}$ be (disjoint) compact Hausdorff spaces such that there are homeomorphisms $\sigma_i ^j $ from $X_i$ to $X_j$  satisfying $\sigma_j ^k \circ \sigma_i^j = \sigma_i^k$, $\sigma_i^i= \id_{X_i} $;  for  instance consider  $X_i=[i-1, i-2/3]$ and $\sigma_i^j(t)= t+(j-i)$ where $ 1 \le i, j \le n$.   

Then let $X= \cup_{i=1}^n X_i $ and  consider $A=C(X)$ and for $f \in A$ we write $f=(f_1, f_2, \dots, f_n)$  of which each $f_i$ is the element of $C(X_i)$ and set $$P=\{(f_1, f_2, \dots f_n) \in A \mid f_i (\sigma_j^i (t))= f_j(t)  \quad \text{for all} \,\,  i,j \}. $$  
Then we define a conditional expectation from A to P as follows; given  $f=(f_1, \dots, f_n) \in A$, 
\begin{equation} \label{E:CE1}
 E(f) (t)= \frac{\sum_k f_k(\sigma_i^k(t)) }{n} \quad \text{for}\quad t \in X_i.
 \end{equation}
Indeed, if we write $E(f)=(g_1, \dots, g_n)$, then  for $t \in X_j$ 

\[ 
\begin{split}
g_i (\sigma_j^i(t)) &= \frac{\sum_k  f_k ((\sigma_i^k \circ \sigma_j^i )(t))}{n}\\
&=  \frac{\sum_k  f_k (\sigma_j^k (t))}{n}\\
&= g_j(t), \quad 
\end{split}
\]
so it is well-defined. 
Then  define $\displaystyle e_j = \begin{cases}
&1 \quad \text{if $x\in X_j$}\\
&0 \quad \text{otherwise}
\end{cases}$, or $e_j = (0, 0, \dots, 1, 0, \dots, 0)$. 
It follows that  $E(e_j)=\dfrac{1}{n}(1, 1, \dots, 1)$ or $E(e_j)= \dfrac{1}{n} 1_{A}$ for all $j$. 

Let us show that  $E$ is of index-finite type and has the Rokhlin property. 
 
Consider $\{(u_i=ne_i, v_i=e_i) \mid n=1, \dots, n\}$, then it is a quasi-basis for $E$; for $f= (f_1, \dots, f_n)$ 
\[
\begin{split}
\sum_{i=1}^n\ u_i E(v_i f) &= \sum_{i=1}^n ne_i \left( \frac{f_i \circ \sigma_i^1}{n},\dots, \frac{f_i}{n}, \dots, \frac{f_i\circ \sigma_n^i}{n} \right) \\
&=\sum_{i=1}^n  \left(0, \dots, 0, f_i, 0, \dots, 0 \right)=(f_1, \dots, f_n) \\
&=f.
\end{split}
\]
Similarly, $\sum_{i=1}^n E(f u_i ) v_i =f$ is obtained.  Note that $\Index E = \sum_{i=1}^n u_i v_i=n I_A$.  To show that $E$ has the Rokhlin property, we consider 
 $e$ an element of $A^{\infty} $ defined by the following periodic sequence $(e_1, e_2, \dots, e_n, e_1, e_2, \dots, e_n, \dots)$. Since $A$ is commutative, automatically $e\in A_{\infty}\cap A' $. Moreover, 
$$
E^{\infty}(e)= [(E(e_1), E(e_2), \dots, E(e_n),  E(e_1), \dots)]= \frac{1}{n} I_{A_{\infty}}. 
$$ 
Thus, 
$\displaystyle (\Index E) E^{\infty}(e)=I_{A_{\infty}}$. 
Finally, we show that the map 
$\displaystyle A \ni f \mapsto fe \in A_{\infty}$ is injective.  
For $f=(f_1, \dots, f_n)$, 
$fe = [(\mathbf{f}_1, \mathbf{f}_2, \dots,  \mathbf{f}_n, \mathbf{f}_1, \dots )]$ where $\mathbf{f_j}$ is the function such that  $f_j$ on $X_j$ and vanishes elsewhere. Write $fe$ as $[(\mathbf{f}_k)] $ where $\mathbf{f}_{nq+l}= \mathbf{f}_l$ for all $q, l \in \mathbf{N}$. Let $fe=0$, then for any $\epsilon>0$, there is a natural number $N$ such that  $\| \mathbf{f}_k \| < \epsilon$ for all $k \ge N$. Write $N= nq+r$ for some natural numbers $q\ge 0$ and $r < n$. Then we have 
\[ \| \mathbf{f}_{n(q+1)+j}  \| < \epsilon \] 
for $j=1,2, \dots, n$. Consequently, 
\[ \| \mathbf{f}_{j}  \| < \epsilon \] 
for $j=1,2, \dots, n$  which implies $\|f\| < \epsilon$ in $A$. Thus $f=0$ in $A$.
So we are done.   
\end{ex}

%%%%%%%%%%%%%%%%%%%%%%%%%%%%%%%%%%%%%%%%%%%%%%%%%%%%%%%%%%%%%
Note that $C(X_i)$'s  are isomorphic to each other in Example \ref{Ex:1}. Thus if we consider $n$-copies of the same algebra and forget the homeomorphisms between them,  we  may expect to have an inclusion with the Rokhlin property irrelevant of group actions.     

 \begin{ex}\label{Ex:2}
Let $A$ be a unital $C\sp*$-algebra and for any $n\in \mathbb{N}$ with $n\ge 2$  consider $B=A^n= A\oplus A \oplus \cdots \oplus A$,  the  $n$-direct sum of $A$'s. Let $D=\{(x,x,\dots, x)\mid x\in A \}$ be a unital sub-$C\sp*$-algebra of $B$.  We define a conditional expectation $E:B \to D$ by 

 \begin{equation}\label{E:CE2}
      E(x_1, x_2, \cdots, x_n)= \left( \frac{\sum_{i=1}^n x_i}{n}, \cdots,  \frac{\sum_{i=1}^n x_i}{n} \right)
 \end{equation}
  
Let $e_1=(1, 0, \dots, 0), e_2=(0, 1, 0,\dots, 0), \dots, e_n=(0,\dots, 0,1)$ are the canonical central projections for $A^n$, then define $u_i=ne_i$ and $v_i=e_i$ for $i \le i \le n$. We claim that $\{(u_i, v_i) \mid i=1, \dots, n\}$
is a quasi-basis for $E$ as follows; for any $x=(x_i) \in B$, 
\[ \begin{split}
\sum_{i=1}^n\ u_i E(v_ix) &= \sum_{i=1}^n ne_i \left( \frac{x_i}{n}, \frac{x_i}{n}, \dots, \frac{x_i}{n}\right) \\
&=\sum_{i=1}^n n \left(0, \dots, 0, \frac{x_i}{n}, 0, \dots, 0 \right)=(x_1, \dots, x_n)=x.
\end{split}
\]
Similarly,  $\sum_i^n E(xu_i)v_i=x$. Note that $\Index E = \sum_i u_iv_i=n$. 

We claim that $E$ has the Rokhlin property: Note that $E(e_i)= \dfrac{1}{n} 1_B$ 
for every $i=1,\dots, n$.  
We define a Rokhlin projection $e$ by a periodic infinite sequence 
$(e_1, e_2, \dots, e_n, e_1, e_2, \dots, e_n, e_1, \dots)$. 
Since $e_i$'s are central projections, it follows that $e \in B_{\infty}\cap B'$ and 
$E^{\infty}(e)=[(E_(e_j)_{j=1}^{\infty}]=\dfrac{1}{n}1_{B_{\infty}}$, thus 
$(\Index E) E^{\infty}(e)=1_{B_{\infty}}$.  
Finally, we check that the map 
$B \ni x  \mapsto xe \in B_{\infty}$ is injective. 
Let $x=(x_1, x_2, \dots, x_n)$ and $\mathbf{x}_i=(0,\dots, 0, x_i, 0, \dots, 0)=xe_i$ 
for each $i$, then $xe$ is also represented by a periodic sequence 
$(\mathbf{x}_1, \mathbf{x}_2, \cdots, \mathbf{x}_n, \mathbf{x}_1, \dots)$. 
Now suppose that $xe=0$, then for any $\epsilon>0$, 
there is a natural number $N$ such that  
$\| \mathbf{x}_k \| < \epsilon$ for all $k \ge N$. 
Write $N= nq+r$ for some natural numbers $q\ge 0$ and $r < n$. Then we have 
\[ \| \mathbf{x}_{n(q+1)+j}  \| < \epsilon \] 
for $j=1,2, \dots, n$. Consequently, 
\[ \| \mathbf{x}_{j}  \| < \epsilon \] for $j=1,2, \dots, n$  which implies $\|x\| < \epsilon$ in $B$. Thus $x=0$ in $B$.  
\end{ex}

\begin{rem}
If we consider the mod $n$-shift action or the circular action on $B=A^n $ in the above, then we have a $\mathbb{Z}/ n \mathbb{Z}$-action $\alpha$ on $B$ and it is easily checked that $D= B^{\alpha}$ and $E=E_{\alpha}$.  
\end{rem}
The presence of a natural group symmetry- a cyclic permutation of the summands- suggests that integer-index inclusions may inherently admit group-theoretic interpretations. This further supports the idea that integer indices align with classical symmetries. The examples that we have seen so far are canonical in some sense, now we show that there is a nontrivial example from the hyperfinite $\rm II_1$ factor that appeared in Vaughan Jones' celebrated work \cite{Jones}. 

\begin{ex}
Consider $e_1, e_2, \dots$ be a sequence of Jones' projections in a hyperfinite factor $M$ of type $\rm II_1$ which satisfies the following relations; 
\begin{equation}\label{E:relations}
\begin{split}
e_i e_{i \pm 1}e_i&= \tau e_i \\
e_i e_j &= e_j e_i \quad \text{for} \quad | i - j | \ge 2. 
\end{split}
\end{equation}
and  $M= \{e_1, e_2, \dots\}^{''}$, where $\tau^{-1}$ is a fixed constant in $\{ 4\cos ^2 \frac{\pi}{n} \mid n= 3, 4, \dots \}$. In this case, when $\mathcal{B}$ is the $C\sp*$-algebra generated by $1, e_1, e_2, \dots$, and $\mathcal{A}$ is the sub-$C\sp*$-algebra generated by $1, e_2, e_3, \dots$, there exists a conditional expectation $E:\mathcal{B} \to \mathcal{A}$ such that 
\begin{equation}\label{E:CE}
E(xe_1y)=\tau xy \quad \text{for} \quad x, y \in \mathcal{A}.
\end{equation}
Now we consider the $C\sp*$-algebra generated by $1, e_1, e_3, e_4, \dots$, denoted by $B$, and the $C\sp*$-algebra  generated by $1, e_3, e_4, \dots $, denoted by $A$. Then $B \subset \mathcal{B}$ and $A\subset \mathcal{A}$ and $E(B)=A$. Thus if we restrict $E$ on $B$, then we obtain the conditional expectation from $B$ onto $A$ and still denote it by $E$. 

Let $ \displaystyle u_1= \frac{e_1}{\tau}, v_1=e_1, u_2=\frac{1-e_1}{1-\tau}, v_2=1-e_1$. Then $\{(u_1, v_1), (u_2, v_2)\}$ is a quasi-basis for $E$ as follows; for any $x\in A$ 
 \[ \begin{split}
 u_1E(v_1x)+u_2E(v_2 x) &= \frac{e_1}{\tau} E(e_1x) + \frac{1-e_1}{1-\tau}E((1-e_1)x) \\
                                          &= \frac{e_1}{\tau} \tau x + \frac{1-e_1}{1-\tau}(1-\tau)x \quad \text{by (\ref{E:CE})}\\ 
                                          &= e_1 x + (1-e_1)x =x.
                                          \end{split}\]
Similarly, $E(xu_1)v_1+E(xu_2)v_2= x $. Thus $\displaystyle \Index E = \frac{e_1}{\tau}+ \frac{1-e_1}{1-\tau}= \frac{1}{\tau}$ if $\tau= 1-\tau$ or $\tau=1/2$ which is indeed the case $\tau^{-1}= 4\cos^2\frac{\pi}{n} $ and  $n=4$. One more benefit of such a choice is that $E(e_1)= E(1-e_1)=\tau$. 

Now we claim that $E$ has the Rokhlin property;  we take  $e$ in $B^{\infty}$ which is represented by a periodic sequence $(e_1, 1-e_1, e_1, 1-e_1, \dots )$. Then from (\ref{E:relations}) $e\in B^{\infty}\cap B'$. It also follows that $E^{\infty}(e)= \tau I_{B^{\infty}}$. Moreover, as we have seen in the above, $(\Index E) E^{\infty}(e)=I_{B^{\infty}}$.  By the same argument in  Example \ref{Ex:2}, the map $B \ni x \mapsto xe$ is injective. So we are done.   
\end{ex}
\begin{rem}
If we define an automorphism $\alpha$ on $B$ by sending $e_1$ to $1-e_1$ and other elements fixed, then we can easily see that $\alpha$ is of order $2$ and $E=E_{\alpha}$ in this case.  
\end{rem}
This example, derived from subfactor theory, demonstrates that even in  infinite-dimensional settings inclusions with integer indices can often be traced back to a hidden symmetry structure. The choice of conditional expectation in this setting aligns with known results in Jones' theory, where finite-index inclusions frequently correspond to group-like symmetries (e.g., those  arising from crossed products). This provides further evidence that inclusions with integer indices exhibit intrinsic group-theoretic properties, even when not explicitly defined by a group action. 

%%%%%%%%%%%%%%%%%%%%%%%%%%%%%%%%%%%%%%%%%%%%%%%%%%%%%%

\section{Rokhlin inclusions with Non-integer Index}\label{S:4}
In contrast to the examples in the previous section, here we investigate inclusions of unital $C\sp*$-algebras that possess the Rokhlin property but whose  \emph{Watatani indices are not  necessarily integers} and \emph{do not arise from finite group actions}.  These examples stem from genuine quantum symmetries, such as those implemented by injective $\rm II_1$-subfactors or unitary fusion categories, and highlight the existence of Rokhlin-type phenomena beyond classical group symmetry.     

\subsection{ Inclusion from Cuntz algebra $\mc{O}_2$}

\begin{lem}\label{L:Rokhlin2}
Let $A$ be a unital separable simple $C\sp*$-algebra and $P\subset A$ be an inclusion of unital $C\sp*$-algebras and $E: A \to P$ be a conditional expectation of index-finite type. Then $E$ has the Rokhlin property if for every $\epsilon>0$ and a finite subset $\mc{F} \subset A$ there exists a projection  $e$ such that 
\begin{enumerate}
\item $\| ea -ae \| < \epsilon $ for all  $a\in \mc{F}$,
\item $\| 1-(\Index E) E(e) \| < \epsilon$.  
\end{enumerate} 
\end{lem}

Based on a study of fusion rule in sector theory, or fusion category action on $\mc{O}_2$ in modern language, Izumi constructed the following conditional expectation on $\mc{O}_2$ with a finite  index as an irrational number. 
  
\begin{thm} \label{T:Izumi}
Let $\mathcal{O}_2$ be the Cuntz algebra whose generators are $S_1, S_2$ and $\rho: \mc{O}_2 \to \mc{O}_2$ the endomorphism defined by 
\[
\begin{split}
\rho(S_1)&= \frac{1}{d} S_1 + \frac{1}{\sqrt{d}} S_2S_2,\\
\rho(S_2)&= \left( \frac{1}{\sqrt{d}}S_1 - \frac{1}{d} S_2 S_2 \right) S^*_2 +S_2 S_1 S^*_1
\end{split} 
\] where $d$ is the positive real number satisfying $\displaystyle 1= \frac{1}{d}+ \frac{1}{d^2}$. Then the conditional expectation $E_{\rho}: \mc{O}_2 \to \Image \rho$ defined by $E_{\rho}(x)=\rho(S^*_1 \rho(x)S_1)$ has $Index E_{\rho}=d^2$ and its quasi-basis is $\{dS_1^*, dS_1\}$. Moreover, $e=S_1S_1^*$ satisfies that $\Index (E_{\rho}) E_{\rho}(e)=1$. 
\end{thm}
\begin{proof}
See \cite[Example 3.2]{Izumi}. 
\end{proof}

We need Elliott's theorem $\mc{O}_2 \otimes \mc{O}_2 \cong  \mc{O}_2$, which we utilize to construct an inductive system for $\otimes_{i=1}^{\infty} \mc{O}_2$. The following is well-known but the details of the proof is important to us since we repeatedly use this machinary.  

\begin{prop}\label{P:O_2}
Let $\phi: \mc{O}_2 \to \mc{O}_2 \otimes \mc{O}_2$ be the injective $*$-homomorphism given by the second factor embedding or $\phi(x)=1\otimes x$. Then there is a sequence of unitaries $\{ v_n\}_{n=1}^{\infty}$ in $\mc{O}_2\otimes \mc{O}_2$ such that 
\[ \lim_{n\to \infty} \|v_n \phi(a) -\phi(a)v_n \|=0, \quad  \lim_{n\to \infty} \dist (v_n^*b v_n, \phi( \mc{O}_2) ) =0\]
for all $a\in \mc{O}_2$ and all $b \in \mc{O}_2\otimes \mc{O}_2$. Consequently, there is an isomorphism $\psi: \mc{O}_2 \to \mc{O}_2\otimes \mc{O}_2$ which is approximately unitarily equivalent to $\phi$ denoted by $\phi \approx_{u} \psi$.
\end{prop}  
\begin{proof}
$v_n$'s in the above essentially come from an asymptotic central sequence of endomorphisms, $\{\rho_n\}_{n=1}^{\infty}$ (see \cite[Lemma 5.2.3]{Ro}). In fact, for an arbitrary $\delta>0$ there is a unitary $w$ in $\mc{O}_2 \otimes \mc{O}_2$ such that $\| w(1 \otimes S_j )w^*- S_j \otimes 1 \| < \delta $ and  $\| w(1 \otimes S_j^* )w^*- S_j^* \otimes 1 \| < \delta$ for $ j=1,2$ by the uniqueness result for $\mc{O}_2$ or the approximately inner half flip property of $\mc{O}_2$. So we take $w_n= (\id \otimes \rho_n) (w)$. 
Then  
\begin{align*}
\| w_n (1\otimes S_j) - (1\otimes S_j)w_n\| & \to 0 \quad \text{as} \quad  n\to \infty \\
\|w_n (1\otimes S_j^*) - (1\otimes S_j^*)w_n\| & \to 0 \quad \text{as} \quad  n\to \infty 
\end{align*}
for $ j=1,2$. Also, 
\[
\begin{split}
\dist (w_n^* (S_j \otimes 1) w_n, 1\otimes \mc{O}_2) &\le \| w_n^* (S_j \otimes 1) w_n- 1\otimes \rho_n(S_j)\| \\
& \le \| (\id \otimes \rho_n)(w^* (S_j \otimes 1) w)- (\id \otimes \rho_n)(1\otimes S_j) \| \\
&< \delta 
\end{split}
\]
for $j=1,2$. Then the conclusion follows from \cite[Proposition 2.3.5]{Ro}.  We note that $\psi$ is the limit of the maps of the form $w_{n_1}w_{n_2} \cdots w_{n_{k}}\phi(\cdot)w_{n_k}^*\cdots w_{n_2}^*w_{n_1}^* $ for a suitable choice of $\{w_{n_k}\}$ which is a subsequence of $\{w_n\}$.  
\end{proof} 

\begin{cor}\label{C:inductivelimit}
For each $n\ge 2$ we inductively let $\phi_n: \otimes_{i=1}^n \mc{O}_2 \to \mc{O}_2 \otimes (\otimes_{i=1}^{n} \mc{O}_2) $ be the second factor embedding as in Proposition \ref{P:O_2}. Viewing $\otimes_{i=1}^n \mc{O}_2 \cong \mc{O}_2$ via $\psi_{n-1}$, we can apply  Proposition \ref{P:O_2} to the injective map $\phi_n \circ \psi_{n-1}: \mc{O}_2 \to \otimes_{i=1}^{n+1} \mc{O}_2$.  Therefore there is an isomorphism $\psi_n: \mc{O}_2 \to \otimes_{i=1}^{n+1} \mc{O}_2$ which is approximately unitarily equivalent to $\phi_n \circ \psi_{n-1}$. 
\end{cor}
 
\begin{rem}\label{R:inductivelimit}
In the above inductive process, in each $n$-th step,  it is identical to use the approximately inner half flip property of $\mc{O}_2$ and an asymptotically central sequence of endomorphisms $\rho_n$'s  while replacing $S_j$ and $S_j^*$ by $\psi_{n-1}(S_j)$ and $\psi_{n-1}(S_j^*)$ respectively.  
\end{rem}

We put $\psi_1=\psi$ corresponding to $\phi_1= \phi$ in Proposition \ref{P:O_2}.  If we inductively replace $\psi_n $ by $\Ad W_n \circ \psi_n$ by suitable unitary $W_n$ in $\otimes_{i=1}^{n+1} \mc{O}_2$,  we have the following diagram for the inductive limit of a system $(\otimes_{i=1}^n \mc{O}_2, \phi_n: \otimes_{i=1}^n \mc{O}_2 \to \otimes_{i=1}^{n+1} \mc{O}_2)$ which is approximate intertwining in the sense of Elliott \cite{Ell}(here we still abuse $\psi_n$ for $\Ad W_n \circ \psi_n$ without confusion); 
\bigskip
\begin{equation}\label{E:approximateintertwining}
\xymatrix{
& \mc{O}_2 \ar@{->}[rr]^{\phi_1 } \ar@<1ex>[dr]^{\id}
&& \mc{O}_2\otimes \mc{O}_2  \ar@{->}[rr]^{\phi_2 }  \ar@{->}[dr]^{\psi_1^{-1}}&& \mc{O}_2\otimes \mc{O}_2\otimes \mc{O}_2 \ar@{->}[r] & \cdots \ar@{->}[r] & \otimes_{i=1}^{\infty} \mc{O}_2 \ar@<0.1ex> @{-->}[d]\\ 
\mc{O}_2\ar@{=}[ur]^{\id}
\ar@{=}[rr]^{id} && \mc{O}_2 \ar@{->}[ur]^{\psi_1}\ar@{=}[rr]^{id} && \mc{O}_2 \ar@{->}[ur]^{\psi_2}  \ar@{=}[rr]^{id}  &  & \cdots \ar@{=}[r]^{id}& \mc{O}_2  \ar@<1ex> @{-->}[u]^{\psi}}
\end{equation}

%\begin{equation} \begin{CD}
%O_2 @>\phi_1 >> O_2 \otimes O_2 @>>> \dots @>>> \otimes_{i=1}^n O_2 @> \phi_n>> \otimes_{i=1}^{n+1}O_2 @>>> \dots @>>> \otimes_{i=1}^{\infty}O_2\\
%@|                   @V \psi_1^{-1} VV    @.    @V\psi_{n-1}^{-1}VV    @V\psi_{n}^{-1}VV      @.          @VVV \\
%O_2 @= O_2  @= \dots @= O_2 @= O_2 @= \dots @= O_2
%\end{CD} 
%\end{equation}
\bigskip
\begin{prop}\label{P:CE}
Consider $E_{\rho}(x)=\rho(S_1^* \rho(x)S_1)$ where $\rho$ is  the endomorphism of $ \mc{O}_2$
given in Theorem \ref{T:Izumi}.  Then we define $E_{n+1}: \otimes_{i=1}^{n+1} \mc{O}_2 \to \psi_n(\Image \rho)$ by $E_{n+1}= \psi_{n} \circ E_{\rho} \circ \psi_n^{-1}$ for each $n \in \mathbb{N}$ and $E: \otimes_{i=1}^{\infty} \mc{O}_2 \to \mc{O}_2$ by $\lim_{n\to \infty}E_n $.  It follows that $\Index E =\Index E_{\rho}$ and $E$ has the Rokhlin property. 
\end{prop}
\begin{proof}
Note that for each $n$ $E_{n+1}$ is the conditional expectation such that $\Index E_{n+1} =\Index E_{\rho}$. Since we can easily check that $E_{n+1} \circ \phi_n = \phi_{n}\circ E_{n}$ approximately, $E$  is well-defined as the inductive limit of $E_n$ (see \cite[p.181 2.3]{Ell}). It remains to show that $E$ has the Rokhlin property: for a finite set $F \subset \otimes_{i=1}^{\infty} \mc{O}_2$ and $\epsilon >0$, there is a $M \in \mathbb{N}$ such that $F\subset  \otimes_{i=1}^{M} \mc{O}_2$.
Then we consider a projection $\widehat{e}=\overbrace{1 \otimes \cdots \otimes1}^{M \rm\ times }\otimes \psi_N (e) \in \otimes_{i=1}^{N+M+1} \mc{O}_2$ for some $N$ where $e=S_1S_1^*$, then $\widehat{e}$ commutes with elements in $F$.  Note that $\widehat{e}=(\phi_{N+M} \circ \phi_{N+M-1}\circ \cdots \circ \phi_{N+1})(\psi_N(e))$. Thus if we choose $N$ large enough,  $\widehat{e}$ is close to $ \psi_{N+M}(e)$  up to $\epsilon$ by the approximate intertwining of the diagram (\ref{E:approximateintertwining}). If necessary, we can choose $N$ such that $\| (\Index E) E (\psi_{N+M}(e))-(\Index E)E_{N+M+1}(\psi_{N+M}(e))\| < \epsilon$.

Then 
\[
\begin{split}
 \| 1-( \Index E) E(\widehat{e})   \| &\le \| 1- (\Index E) E_{N+M+1} (\psi_{N+M}(e))\| \\
 &+\| (\Index E) E (\psi_{N+M}(e))-(\Index E)E_{N+M+1}(\psi_{N+M}(e))\|\\
  &+ \| (\Index E) E(\psi_{N+M}(e)) - (\Index E) E(\widehat{e})\| \\
 &< 0+ \epsilon + (\Index E) \epsilon= (1+\Index E) \epsilon.  
 \end{split}\]
 
Since $\otimes_{i=1}^{\infty} \mc{O}_2 (\cong \mc{O}_2)$ is simple,  we are done by Lemma \ref{L:Rokhlin2}. 
\end{proof}

We need to make sure that the above construction of $P=\Image E \subset \otimes_{i=1}^{\infty} \mc{O}_2$ with the conditional expectation $E$ as in Proposition \ref{P:CE} is irreducible. To show it, we begin with a simple observation.   

\begin{lem}\label{L:irreducible}
Let $P \subset A$ be an inclusion of simple $C^*$-algebras and 
$E: A  \to P$ a conditional expectation from A onto P with $\Index E< \infty$. 
If $\Index E$ is less than 4, then $P \subset A$ is  irreducible, i.e., $P' \cap A \simeq \C$.
\end{lem} 
\begin{proof}
Since $\Index E$ is an element of $A' \cap A \simeq \C$, we have 
$$
\Index E = 4\cos^2\frac{\pi}{n}, \quad n= 3, 4, 5, \cdots
$$
 by \cite[Theorem 2.3.6]{W:index}. Therefore the principal graph of $P \subset A$ is in 
$A_n, D_n, E_n$ by \cite[Theorem 1.4.3]{GHJ}.
 In any cases, there is only one edge connected to $*$ and hence $P' \cap A \simeq \C$
, i.e., $P \subset A$ is irreducible.
\end{proof}

\begin{cor} \label{C:irreducible}
Let $P \subset \otimes_{i=1}^\infty \mc{O}_2$ be an inclusion given in Proposition \ref{P:CE}. Then, $P' \cap \otimes_{i=1}^\infty \mc{O}_2 = \C$. In other words, this inclusion is irreducible.
\end{cor}
\begin{proof}
Let $E:\otimes_{i=1}^\infty \mc{O}_2 \rightarrow P$. Then, 
since $\mathrm{Index} E < 4$, from Lemma~\ref{L:irreducible} we know that $P' \cap \otimes_{i=1}^\infty \mc{O}_2 = \C$.
\end{proof}

\begin{rem} \label{R:irreducible}
Using the irreducibility of $\Image \rho \subset \mc{O}_2$ in Theorem \ref{T:Izumi} we can deduce the irreducibility of $P \subset \otimes_{i=1}^\infty \mc{O}_2$ directly; note that $P= \psi(\mc{O}_2)$ where $\psi$ is the isomorphism obtained in the diagram (\ref{E:approximateintertwining}).  
Let $y \in P' \cap \otimes_{i=1}^\infty \mc{O}_2$. Then, there is a $x \in \mc{O}_2$ such that $y = \psi(x)$. Since $y \psi(\rho(a)) = \psi(\rho(a))y$ for any $a \in \mc{O}_2$, $\psi(x)\psi(\rho(a)) = \psi(\rho(a))\psi(x)$, that is, $x\rho(a) = \rho(a)x$ for any $a \in \mc{O}_2$. Since $\rho( \mc{O}_2)'\cap \mc{O}_2 = \C$ \cite[Theorem~5.6(1)]{Izumi}, $x \in \C$. Hence $y = \psi(x) \in \C$. So, $P' \cap \otimes_{i=1}^\infty \mc{O}_2 = \C$.
\end{rem}

\bigskip 
\subsection{ Inclusion of Injective $\rm II_1$-factors}

\begin{defn}
Let $M$ be a von Neumann algebra with a unique tracial state $\tau$ and $\omega$ be a free ultra filter on $\N$. 
Let $\Pi_{n\in M} M= \{(x_n) x_n \in M, \sup_{n\in\N}\|x_n\| < \infty\}$ and  $N_\omega^{(2)} = \{(x_n) \in \Pi_{n\in \N}M | \lim_{n\rightarrow \omega}\|x_n\|_{2,\tau} = 0\}$. Then, define $M^\omega = \Pi_{n\in \N}M/N_\omega^{(2)}$, where $\|x\|_{2,\tau} = \tau(x^*x)^\frac{1}{2}$.
\end{defn}

\vskip 2mm

\begin{rem}
Note that $M^\omega$ becomes a von Neumann algebra (\cite[Lemma A.9]{BO}). Moreover, if $M$ is factor, then, $M^\omega$ is also factor(\cite[Excercise~A5]{BO}).
\end{rem}

When $N \subset M$ is given where $(M, \tau)$ is a finite von Neumann algebra, recall that there is a unique canonical trace preserving conditional expectation $E_{\tau} : M \to N$ that is the orthogonal projection of $L^2(M, \tau)$ onto the closed subspace $L^2(N, \tau)$. 

\bigskip

From the observation in \cite[Theorem 5.3]{ELPW 2010} the following is the suitable definition
of the Rokhlin property for inclusions of factor-subfactor. 

\begin{defn}
 Let $N \subset M$ be of irreducible injective $\rm II_1$-factors of index-finite
type and let $\omega$ be a free ultrafilter. Then, we say that an inclusion $N \subset M$ has the
Rokhlin property if there exists a projection $e \in  M^{\omega}\cap M'$ such that $[M:N] E^{\omega}_{\tau}(e) = 1
$ where $\Index E_{\tau}=[M:N]$.
\end{defn}

\begin{rem}
 The above definition is independent of an ultrafilter $\omega$, that is, if
there exists a projection $e \in M^{\omega}\cap M'$  such that $[M:N]E^{\omega}(e) = 1$, then, for another
ultrafilter $\omega'$ there exists a projection $ f  \in M^{\omega'}\cap M'$
 such that $[M:N]E^{\omega'} (f) = 1$.
\end{rem}

Let us consider $N\subset M$ of irreducible $\rm II_1$- factors as a (quantum) symmetry in operator algebras which produces  non-integer indices.    
 \begin{thm}\label{T:rokhlin}
Let $N \subset M$ be of irreducible injective $\rm II_1$-factors of index-finite type.
Then %for any ultrafileter $\omega$ %$N_\tau^\omega \subset M_\tau^\omega$
$N \subset M$
has the Rokhlin property. That is, 
there exists a projection $e \in M^{\omega} \cap M'$ such that 
$[M:N] E^{\omega}_\tau(e) = 1$. 
\end{thm}

The following observation is useful to prove Theorem \ref{T:rokhlin} which provides the analogous toolkit for the inductive limit with respect to tracial 2-norm $\| \cdot \|_{2, \tau}$.

\vskip 2mm

\begin{defn}
A C*-subalgebra $A$ of a C*-algebra $B$ is said to be relative LF in $B$ if for every finite subset $F \subset  A$ and every $\epsilon > 0$ there is a finite-dimensional $C^*$-algebra $D$ of $B$ such that for every $a \in F$, $\mathrm{dist}(a, D) < \epsilon$.
\end{defn}

\vskip 2mm

\begin{thm}\label{thm:LHL}\cite[Theorem 4]{LHL 2023}
Let $A$ be a unital $C\sp*$-algebra and suppose that $A$ is relative LF in $A^{**}$ and $M$ is a finite von Neumann algebra. If $\rho, \pi:A \rightarrow M$ are unital *-homomorphisms, then the following are equivalent:
\begin{enumerate}
\item
$\pi$ and $\rho$ is approximately unitarily equivalent,
\item
$M-\mathrm{rank}(\pi(a)) = M-\mathrm{rank}(\rho(a))$ for every $a \in \mathcal{A}$
\item
$\Phi_{M} \circ \pi = \Phi_{M} \circ \rho$
\end{enumerate}
where $M-\mathrm{rank}(T)$ is the Murray von Neumann equivalent class of projection $R(T)$ onto the closure of the range of $T$ and $\Phi_{M}$ is a faithful normal tracial conditional expectation from $M$ to its center $\mathcal{Z}(M)$.
\end{thm}

\vskip 2mm

\begin{rem}
\begin{enumerate}
\item
If $M$ is a factor and $A$ has a unique tracial state, any two *-homomorphism
$\pi, \rho$ are unitarily equivalent from the condition in Theorem~\ref{thm:LHL} (3).
\item
If $A$ is an ASH algebra, $A$ is relative LF in $A^{**}$ \cite{LHL 2023}.
\end{enumerate}
\end{rem}

\vskip 2mm

\begin{cor}
Let $N$ be an injective $\rm{II}_1$ factor and $M$ be a factor with a tracial state $\tau$.
Then, any two *-homomorphisms $\pi, \rho: N \rightarrow M$ are unitarily equivalent in the weak *-topology.
That is, there is a sequence $\{u_n\}_{n\in \N}$ of unitary operators in $M$ such that 
$$
\|u_n\pi(a)u_n^*-\rho(a)\|_{2,\tau} = \tau((u_n\pi(a)u_n*-\rho(a))^*(u_n\pi(a)u_n^* - \rho(a)))^{1/2} \rightarrow 0$$
as $n \rightarrow \omega$.
\end{cor}

\vskip 2mm

%\vskip 2mm

The following might be well-known to experts. 

\vskip 2mm

\begin{lem}\label{lem:roklin projection}\cite[Proposition~1.9]{PP 1986}
Let $N \subset M$ be of irreducible injective $\rm II_1$-factor of index-finite type.
Then, %$N_\tau^\omega \subset M_\tau^\omega$
%$N \subset M$
$M$ has a projection $e \in M$ such that 
$\displaystyle [M:N]E(e) = I$. 
\end{lem}

\begin{rem}\label{rmk:unitarily equivalent}
By  Connes's classification Theorem \cite[Theorem~1]{AC 1976},  $M \otimes \mathcal{R} \cong M$  for any injective $\rm II_1$-factor  $M$ where $\mathcal{R}$ is the hyperfinite injective $\rm{II}_1$-factor. Thus  there is an isomorphism 
$\psi:\mathcal{R} \rightarrow \mathcal{R} \overline{\otimes} \mathcal{R}$. 
Hence, from Theorem~\ref{thm:LHL} there is an injective *-homomorphism $\phi_n: \overline{\otimes}_{i=1}^n\mathcal{R} \rightarrow \mathcal{R} \overline{\otimes} (\overline{\otimes}_{i=1}^n\mathcal{R})$ by $\phi_n(x) = 1 \otimes x$ such that 
$\phi_n \circ \psi_{n-1}$ is approximately unitarily equivalent to $\psi_n$ in the topology induced by $\|\quad \|_{2,\tau}$, where 
$\psi_{n-1}:\mathcal{R} \rightarrow \overline{\otimes}_{i=1}^{n}\mathcal{R}$ and $\psi_{n}:\mathcal{R} \rightarrow \overline{\otimes}_{i=1}^{n+1}\mathcal{R}$ are isomorphisms.
\end{rem}

\vskip 2mm

{\bf Proof of Theorem~\ref{T:rokhlin}}:

%\color{red}
We may assume that $M$ has a unique tracial state $\tau$, $M = \mathcal{R}$ and $N \cong \mathcal{R}$.
Note that $\Index E_\tau = [M:N]$.

We put $\psi_1=\psi$ corresponding to $\phi_1= \phi$ in Remark \ref{rmk:unitarily equivalent}.  If we inductively replace $\psi_n $ by $\Ad W_n \circ \psi_n$ by suitable unitary $W_n$ in $\overline{\otimes}_{i=1}^{n+1}\mathcal{R}$,  we have the following diagram for the inductive limit of a system $(\overline{\otimes}_{i=1}^n \mathcal{R}, \phi_n: \overline{\otimes}_{i=1}^n \mathcal{R} \to \overline{\otimes}_{i=1}^{n+1} \mathcal{R})$ which is approximate intertwining in the topology induced by $\|\quad \|_{2,\tau}$ in the sense of Elliott (here we still abuse $\psi_n$ for $\Ad W_n \circ \psi_n$ without confusion); 
\bigskip
\begin{equation}\label{E:approximateintertwining}
\xymatrix{
&\mathcal{R} \ar@{->}[rr]^{\phi_1 } \ar@<1ex>[dr]^{\id}
&& \mathcal{R} \overline{\otimes} \mathcal{R}  \ar@{->}[rr]^{\phi_2 }  \ar@{->}[dr]^{\psi_1^{-1}}&& \mathcal{R} \overline{\otimes} \mathcal{R} \overline{\otimes} \mathcal{R} \ar@{->}[r] & \cdots \ar@{->}[r] & \overline{\otimes}_{i=1}^{\infty} \mathcal{R} \ar@<0.1ex> @{-->}[d]\\ 
\mathcal{R} \ar@{=}[ur]^{\id}
\ar@{=}[rr]^{id} && \mathcal{R} \ar@{->}[ur]^{\psi_1}\ar@{=}[rr]^{id} && \mathcal{R} \ar@{->}[ur]^{\psi_2}  \ar@{=}[rr]^{id}  &  & \cdots \ar@{=}[r]^{id}& \mathcal{R}  \ar@<1ex> @{-->}[u]^{\psi}}
\end{equation}

%\begin{equation} \begin{CD}
%O_2 @>\phi_1 >> O_2 \otimes O_2 @>>> \dots @>>> \otimes_{i=1}^n O_2 @> \phi_n>> \otimes_{i=1}^{n+1}O_2 @>>> \dots @>>> \otimes_{i=1}^{\infty}O_2\\
%@|                   @V \psi_1^{-1} VV    @.    @V\psi_{n-1}^{-1}VV    @V\psi_{n}^{-1}VV      @.          @VVV \\
%O_2 @= O_2  @= \dots @= O_2 @= O_2 @= \dots @= O_2
%\end{CD} 
%\end{equation}
%\bigskip
%\begin{thm}\label{P:CE}
%Consider $E_{\rho}(x)=\rho(S_1^* \rho(x)S_1)$ where $\rho$ is  the endomorphism of $O_2$
%given in Theorem \ref{T:Izumi}.  Then we define $E_{n+1}: \otimes_{i=1}^{n+1} O_2 \to \psi_n(\Image \rho)$ by $E_{n+1}= \psi_{n} \circ E_{\rho} \circ \psi_n^{-1}$ for each $n \in \mathbb{N}$ and $E: \otimes_{i=1}^{\infty} O_2 \to O_2$ by $\lim_{n\to \infty}E_n $.  It follows that $\Index E =\Index E_{\rho}$ and $E$ has the Rokhlin property. 
%\end{thm}

%\begin{proof}

Note that for each $n$ $E_{n+1}: \otimes_{n=1}^{n+1}\mathcal{R} \rightarrow \psi_n(N)$ by $E_{n+1} = \psi_n \circ E_\tau \circ \psi_n^{-1}$ for each $n \in \N$ is the conditional expectation such that $\mathrm{Index} E_{n+1} =\mathrm{Index}E_{\tau}$. Since we can easily check that $E_{n+1} \circ \phi_n = \phi_{n}\circ E_{n}$ approximately, $E$  is well-defined as the inductive limit of $E_n$ (see \cite[pp181 2.3]{Ell}). It remains to show that $E$ has the Rokhlin property: for given a finite set $F \subset \overline{\otimes}_{i=1}^{\infty}\mathcal{R}$ and $\epsilon >0$, there is a $p \in \mathbb{N}$ such that $F\subset  \overline{\otimes}_{i=1}^{p} \mathcal{R}$.
Then we consider a projection $\widehat{e}=\overbrace{1 \otimes \cdots \otimes1}^{p \rm\ times }\otimes \psi_N (e) \in \overline{\otimes}_{i=1}^{q+p+1}\mathcal{R}$ for some $q$ where $e$ is a projection such that $E_\tau(e) = \frac{1}{[M:N]}$, then $\widehat{e}$ commutes with elements in $F$.  Note that $\widehat{e}=(\phi_{q+p} \circ \phi_{q+p-1}\circ \cdots \circ \phi_{q+1})(\psi_N(e))$. Thus if we choose $q$ large enough,  $\widehat{e}$ is close to $ \psi_{q+p}(e)$  up to $\epsilon$ in the topology induced by $\|\quad\|_{2,\tau}$ by the approximate intertwining of the diagram (\ref{E:approximateintertwining}). If necessary, we can choose $q$ such that $\| (\mathrm{Index} E_\tau) E (\psi_{q+p}(e))-(\mathrm{Index} E_\tau)E_{q+p+1}(\psi_{q+p}(e))\|_{2,\tau} < \epsilon$. \\
Then 
\[
\begin{split}
 \| 1-( \mathrm{Index} E_\tau) E(\widehat{e}) \|_{2,\tau} &\le \| 1- (\mathrm{Index} E_\tau) E_{q+p+1} (\psi_{q+p}(e))\|_{2,\tau} \\
 &+\| (\mathrm{Index} E_\tau) E (\psi_{q+p}(e))-(\mathrm{Index} E_\tau)E_{q+p+1}(\psi_{q+p}(e))\|_{2,\tau}\\
  &+ \| (\mathrm{Index} E_\tau) E(\psi_{q+p}(e)) - (\mathrm{Index} E_\tau) E(\widehat{e})\|_{2,\tau} \\
 &< 0+ \epsilon + (\mathrm{Index} E_\tau) \epsilon= (1+\mathrm{Index} E_\tau) \epsilon.  
 \end{split}
 \]
 
Since $\overline{\otimes}_{i=1}^{\infty}\mathcal{R} (\cong \mathcal{R})$ is simple,  we are done by Lemma \ref{L:probRokhlin}. 
%\end{proof}

\color{black}%%%%%%%%%%%%%%%%%%%%%%%%%%%%%%%%%%%%%%%%%%%%%%%%%%%%%%%%%%%%%%%%

\section{The tracial Rokhlin inclusions}

Next we turn to the tracial case.  We fix a free ultrafilter $\omega$ on $\mb{N}$, and  set $c_{\omega}(A)= \{ (a_n)_n \mid \lim_{n \to \omega} \|a_n\|=0 \}$.  Then $A^{\omega}=\prod_{i=1}^{\infty} A / c_{\omega}(A)$ is called the ultrapower $C\sp*$-algebra of $A$ with the norm of $a \in A^{\omega}$ is given by $\lim_{n \to \omega} \|a_n\|$, where $(a_n)_n$ is a representing sequence of $a$.  Note that the following definition is originally presented using $A^{\infty}$ but we correct it using the ultrapower as explained in \cite{LO}.  
\begin{defn} \cite[Definition 4.2]{OT1} \label{D:tracialRokhlin}
Let  $P \subset A$ be an inclusion of unital $C\sp*$-algebras and let $E \colon A \to P$ be a conditional expectation of index-finite type. We denote by $E^{\omega}$ the canonical conditional expectation from $A^{\omega}$ to $P^{\omega}$ induced by $E$. A conditional expectation $E$ is said to have the tracial Rokhlin property if for any nonzero positive $z \in A^{\omega}$ there exists a projection  $ e \in A' \cap A^{\omega}$ satisfying that $(\Index E)E^{\omega}(e) = g$ is a projection, and $1-g$ is Murray–von Neumann equivalent to a projection in the hereditary subalgebra of $A^{\omega}$ generated by $z$, and a map  $A \ni x \mapsto xe$ is injective. In this case, $e$ is called a Rokhlin projection.
\end{defn} 

More often, the following condition is more practical to check the tracial Rokhlin property.  
\begin{lem}\label{L:tracialRokhlin}
Let $A$ be a separable simple $C\sp*$-algebra and $P\subset A$ be an inclusion of unital $C\sp*$-algebras and $E: A \to P$ be a conditional expectation of index-finite type. A conditional expectation $E$ has the tracial Rokhlin property if for every $\epsilon>0$, a finite subset $\mc{F} \subset A$, and a positive nonzero element $a \in A$ there exists a projection (positive contraction) $e$ such that  
\begin{enumerate}
\item $ (\Index E) E(e)$ is an idempotent,
\item $\| ea -ae \| < \epsilon $ for every $a\in \mc{F}$, 
\item $ 1-(\Index E) E(e)$ is Murray-von Neumann equivalent to a projection in $\overline{aAa}$.   
\end{enumerate} 
\end{lem}
%%%%%%%%%%%%%%%%%%%%%%%%%%%%%%%%%%%%%%%%%%%%%%%%%%%%%%%%%%%%%%%%
\bigskip
We note that the following remarkable result is due to  Kitamura's study on the action of unitary tensor category on Kirchberg algebras  which inherits and develps Izumi's earlier work on Fibonacci fusion category for Theorem \ref{T:Izumi}. Hence again we exploit genuine quantum symmetries to give rise to inclusions with a real Watatani index. 
\begin{thm}\cite[Proposition 5.9]{Ki} \label{T:CE2}
Let $A$ be a unital Kirchberg algebra and $d \in [4, \infty)$. Then there is a conditional expectation $E$ from $A$ onto the image of an endomorphism $\iota: A \to A$ such that $\iota(A)'\cap A=\mb{C}$ and $\Index E =d$. 
\end{thm}
\begin{proof} 
We emphasize that $E: \mc{K}(\alpha(\pi)) \to A$ of index-finite type  is obtained via the action $(\alpha, u)$  of a  unitary tensor category  on $A$ with a  nonzero object $\pi \in \mc{C}$ and $\alpha(\pi)$ as an $(A, A)$-correspondence (see \cite[Lemma 1.26]{KW}) and $\mc{K}(\alpha(\pi)) \cong A$ by Kirchberg-Phillips theorem.   
\end{proof}

Due to R{\o}rdam \cite{Ro},  it is also known that $\mc{O}_{\infty} \cong \mc{O}_{\infty} \otimes \mc{O}_{\infty}$, and more importantly,  it has approximately half inner flip (see \cite[Section 1]{TW}), and thus we have the following result.  

\begin{prop}\label{P:ILS}
Let $\phi: \mc{O}_{\infty} \to \mc{O}_{\infty} \otimes \mc{O}_{\infty}$ be the injective $*$-homomorphism given by the second factor embedding or $\phi(x)=1\otimes x$. Then there is a sequence of unitaries $\{ v_n\}_{n=1}^{\infty}$ in $\mc{O}_{\infty}\otimes \mc{O}_{\omega}$ such that 
\[ \lim_{n\to \infty} \|v_n \phi(a) -\phi(a)v_n \|=0, \quad  \lim_{n\to \infty} \dist (v_n^*b v_n, \phi( \mc{O}_{\infty}) ) =0\]
for all $a\in \mc{O}_{\omega}$ and all $b \in \mc{O}_{\infty}\otimes \mc{O}_{\infty}$. Consequently, there is an isomorphism $\psi: \mc{O}_{\omega} \to \mc{O}_{\infty}\otimes \mc{O}_{\infty}$ which is approximately unitarily equivalent to $\phi$ denoted by $\phi \approx_u \psi$.
\end{prop} 

Following the same argument in Corollary \ref{C:inductivelimit},  with $\phi_n: \otimes_{i=1}^n \mc{O}_{\infty} \to \otimes_{i=1}^{n+1} \mc{O}_{\infty}$ the second factor imbedding. we  have the following diagram for the inductive limit of a system  which is approximate intertwining in the sense of Elliott. 
\bigskip
\begin{equation}\label{E:approximateintertwining2}
\resizebox{\textwidth}{!} {\xymatrix{
& \mc{O}_{\infty} \ar@{->}[rr]^{\phi_1 } \ar@<1ex>[dr]^{\id}
&& \mc{O}_{\infty}\otimes \mc{O}_{\infty}  \ar@{->}[rr]^{\phi_2 }  \ar@{->}[dr]^{\psi_1^{-1}}&& \mc{O}_{\infty}\otimes \mc{O}_{\infty} \otimes \mc{O}_{\infty} \ar@{->}[r] & \cdots \ar@{->}[r] & \otimes_{i=1}^{\infty} \mc{O}_{\infty} \ar@<0.1ex> @{-->}[d]\\ 
\mc{O}_{\infty}\ar@{=}[ur]^{\id}
\ar@{=}[rr]^{id} && \mc{O}_{\infty} \ar@{->}[ur]^{\psi_1}\ar@{=}[rr]^{id} && \mc{O}_{\infty} \ar@{->}[ur]^{\psi_2}  \ar@{=}[rr]^{id}  &  & \cdots \ar@{=}[r]^{id}& \mc{O}_{\infty}  \ar@<1ex> @{-->}[u]^{\psi}}}
\end{equation}

\begin{lem}\label{L:RokhlinProjection}
Let $A$ be a unital purely infinite simple $C\sp*$-algebra, $X$ nondegenerate $(A, A)$-correspondence or Hilbert $A$-bimodule, and $E: \mc{K}(X) \to A$ a conditional expectation of finite index such that $E(\Theta_{x_1, x_2})= (\Index E)^{-1} \langle x_1, x_2 \rangle$ for $x_1, x_2 \in X$. Then for any projection $p$, there is a projection of the form $\Theta_{y, y}$ such that $(\Index E) E(\Theta_{y, y})=p$.   Consequently,  there is a projection of the form $\Theta_{y', y'}$  such that $1- (\Index E) E(\Theta_{y', y'})=p$. 
\end{lem}
\begin{proof}
Note that if $\langle x, x \rangle$ is idempotent, so is $\Theta_{ z, z }$ where $z=x\langle x , x \rangle$. Therefore we have 
\[(\Index E) E(\Theta_{z, z})= \langle x, x \rangle. \]
Since $A$ is simple and purely infinite, there exists $b \neq 0$ such that $1_A = \langle x \cdot b, x \cdot b \rangle$. Hence for any projection $p$, 
\[ p= \langle x \cdot (bp), x\cdot (bp) \rangle \]
Put $y= x\cdot (bp) \langle x\cdot (bp), x \cdot (bp) \rangle= x \cdot (bp)$, and we obtain  $(\Index E) E(\Theta_{y, y})=p$. 
\end{proof}
\begin{prop}\label{P:TCE}
Consider the conditional expectation $E_{\iota}$ as in Theorem \ref{T:CE2} from $\mc{O}_{\infty}$ onto $\iota(\mc{O}_{\infty})$ whose Watatani index is $d$, an irrational number. We define $E_{n+1}: \otimes_{i=1}^{n+1} \mc{O}_{\infty} \to \psi_n(\Image \iota)$ by $E_{n+1}= \psi_{n} \circ E_{\iota} \circ \psi_n^{-1}$ for each $n \in \mathbb{N}$ and $E: A=\otimes_{i=1}^{\infty} \mc{O}_{\infty} \to   \mc{O}_{\infty}$ by the inductive limit of $E_n$ using diagram (\ref{E:approximateintertwining2}).  Then $\Index E =\Index E_{\iota}$ and $E$ has the tracial Rokhlin property.  
\end{prop}
\begin{proof}
For each $n$, $E_{n+1}$ is the conditional expectation such that $\Index E_{n+1} =\Index E_{\iota}$. Since we can easily check that $E_{n+1} \circ \phi_n = \phi_{n} \circ E_{n}$ approximately, $E$  is well-defined as the inductive limit of $E_n$. It remains to show that $E$ has the tracial Rokhlin property; for given $\epsilon >0$, a finite set $F \subset \otimes_{i=1}^{\infty} \mc{O}_{\infty}$, and a positive element $a \in A$ , we may assume that  there is a large enough number $M \in \mathbb{N}$ such that $F\subset  \otimes_{i=1}^{M} \mc{O}_{\infty}$. We choose  a nonzero projection $p \in \mc{O}_{\infty}$. Note that  $\psi_{M+k}(p) \lesssim a$ in $A$ for any $k \in \mb{N}$ since both are nonzero and $A$ is a Kirchberg algebra.   Then we consider a projection $\widehat{e}=\overbrace{1 \otimes \cdots \otimes1}^{M \rm\ times }\otimes \psi (e) \in \otimes_{i=1}^{\infty} \mc{O}_{\infty}$  where $e=\Theta_{y, y}$ such that $1- (\Index E) E_{\iota}(\Theta_{y,y})=p$ by Lemma \ref{L:RokhlinProjection}. It is obvious that $\widehat{e}$ commutes with elements in $F$.  Note that $\widehat{e}$ is the limit of  $(\phi_{N+M} \circ \phi_{N+M-1}\circ \cdots \circ \phi_{N+1})(\psi_N(e))$. Moreover,  if we choose $N$ large enough,  $\widehat{e}$ is close to $ \psi_{N+M}(e)$  up to $\epsilon$ by the approximate intertwining of the diagram (\ref{E:approximateintertwining}). Hence, 
 it follows that $1-( \Index E) E(\widehat{e})$ is the norm limit of  $(1- (\Index E) E_{k+M+1} (\psi_{k+M}(e)))$ as $k \to \infty$. 
 Indeed, we can choose $N$ such that $\| (\Index E) E (\psi_{N+M}(e))-(\Index E)E_{N+M+1}(\psi_{N+M}(e))\| < \epsilon$. \\
Then 
\[
\begin{split}
 &\| 1-( \Index E) E(\widehat{e}) - (1- (\Index E) E_{N+M+1} (\psi_{N+M}(e))) \| \\
 &\le \| (\Index E) E (\psi_{N+M}(e))-(\Index E)E_{N+M+1}(\psi_{N+M}(e))\|\\
  &+ \|  (\Index E) E(\widehat{e})- (\Index E) E(\psi_{N+M}(e))\ | \\
 &< (1+\Index E) \epsilon.  
 \end{split}\]
 
 But for each $n$,  
 \[ (1- (\Index E) E_{n+M+1} (\psi_{n+M}(e)))=\psi_{n+M}(p) \lesssim a. \]
 
 By \cite[II. 3.4.7]{Bl}, $ 1-( \Index E) E(\widehat{e}) \lesssim a$, so we are done by Lemma \ref{L:tracialRokhlin}.  
\end{proof}

It is worth to observe that the above inclusion is irreducible since the index is $\ge 4$. 
\begin{cor} \label{C:irreducible2}
Let $P \subset \otimes_{i=1}^\infty \mc{O}_{\infty}$ be an inclusion given in Proposition \ref{P:ILS}. Then, $P' \cap \otimes_{i=1}^\infty \mc{O}_{\infty} = \C$. In other words, this inclusion is irreducible.
\end{cor}
\begin{proof}
Using the irreducibility of $\Image \iota \subset \mc{O}_{\infty}$ in Theorem \ref{T:CE2} we can deduce the irreducibility of $P \subset \otimes_{i=1}^\infty \mc{O}_{\infty}$ directly following the same argument in Remark
\ref{R:irreducible}.
\end{proof}

%%%%%%%%%%%%%%%%%%%%%%%%%%%%%%%%%%%%%%%%%%%%%%%%%%%%%%%%%%%%%%%%
 \section{Tracial Rokhlin Type Property}

We first recall the notion of  tracial (topological) rank zero of H. Lin \cite{Lin1, Lin2} which is slightly modified from the original one.    
\begin{defn}\label{D:TR} (\cite[Proposition 2.3]{Phillips}) 
Let  $A$ be a simple separable unital $C\sp*$-algebra.  Then $A$ has  tracial rank zero if  and only if  for every finite set $F \subset\subset A$, every $\epsilon > 0$, any nonzero positive element $x \in A$,  there exist projection $p$ and a finite dimensional subalgebra $E \subset pAp$ (that is, $p$ is the identity of $E$) such that 
\begin{enumerate}
\item
$\|pa -ap \| < \epsilon$ for all $a \in F$,
\item
For every $a \in F$, there exists $b \in E$  such that $\|pap - b\| < \epsilon$,
\item
$1-p$ is  Murray-von Neumann equivalent to a projection in $\overline{xAx}$.  
\end{enumerate}
\end{defn}

\begin{prop} $($\cite[Theorem 2.5]{Phillips}$)$
Let $A$ be a simple separable unital tracial rank zero $C\sp*$-algebra. Then $A$ has real rank zero and stable rank one. Moreover, the order on projections in $A$ is determined by the traces $($see \cite[Definition 2.4]{Phillips}$)$.    
\end{prop}

\begin{lem}$($\cite[Lemma 5.2]{ELPW 2010}$)$
When $A$ is a finite, infinite dimensional, simple, separable, unital C*-algebra with Property (SP) (every nonzero hereditary subalgebra contains a nonzero projection) and such that the order on projections over $A$ is determined by traces, we can rephrase the tracial Rokhlin property  of the finite group action $\alpha: G \curvearrowright A$ as follows:
 $\alpha$ has the tracial Rokhlin property if for every finite set $F \subset A$, every $\epsilon > 0$ there exist projections $\{e_g\}_{g\in G}$ such that 
\begin{enumerate}
\item
$e_g \perp e_h$ when $g \not= h$,
\item
$\|\alpha_g(e_h) - e_{gh}\| < \epsilon$ for all $g, h \in G$,
\item
$\|e_gy - ye_g\| < e$ for all $g \in G$ , $y \in F$,
\item
$\tau(1 - \sum_{g\in G}e_g) < \epsilon$ for all $\tau \in \mathrm{T}(A)$,
\end{enumerate}
where $T(A)$ denotes the tracial state space of $A$. 
\end{lem}

When $A$ has tracial rank zero, we have the following characterization.

\vskip 2mm

\begin{thm}\label{thm:ELPW}$($S. Echterhoff, W. L\"{u}ck, N. C.Phillips, and S. Walters \cite{ELPW 2010}$)$
Let $A$ be an infinite dimensional separable unital $C\sp*$-algebra with tracial rank zero. Let $\alpha:G \rightarrow \mathrm{Aut}(A)$ be an action of a finite group $G$ on $A$. Then $\alpha$ has  the tracial Rokhlin property if and only if for every finite set $S \subset A$ and every $\epsilon > 0$, there exist orthogonal projections $e_g \in A$ for $g \in G$ such that
\begin{enumerate}
\item
$\|\alpha_g(e_h) - e_{gh}\|_{2,\tau} < \epsilon$ for all $g, h \in G$ and all $\tau \in T(A)$.
\item
$\|[e_g, a]\|_{2,\tau} < \epsilon$ for all $g \in G$, all $a \in S$, and all $\tau \in T(A)$.
\item
$\sum_{g \in G}e_g = 1$.
\end{enumerate}
\end{thm}

\vskip 2mm

%\begin{dfn} \label{dfn:tracial rokhlin} (cf. \cite[Definition 4.2]{OT 2018})
% Suppose that $A$ is a finite, infinite dimensional simple separable unital C*-algebra with Property (SP) (every nonzero hereditary %subalgebra contains a nonzero projection) and such that the order on projections over $A$ is detemined by traces.
%For any ultrafileter $\omega$ let $C_\omega(A)  = \{(x_n) \in \Pi_{i=1}^\infty A; \tau(x_n^*x_n)^{1/2} \rightarrow 0 (n\rightarrow \omega) %\hbox{for}\ \tau \in T(A)\}$. Then set $\mathcal{A}^\omega = \Pi_{n=1}^\infty A/C_\omega(A)$.

%Let $P \subset A$ an inclusion of unital C*-algebras such that a conditional expectation $E\colon A \rightarrow P$ has a finite index and %$E^\omega$ be the induced map from $A^\omega$ to $P^\omega$.
%It is said to that $E$ has the tracial Rokhlin property if there is a projection $e \in A^\omega \cap A'$ such that $\mathrm{Index}E %E^\omega(e) = 1$.

%if for every nonzero positive element $z \in A_\omega$ there is a projection $e \in A_\omega \cap A'$ such that 
%\begin{enumerate}
%\item
%$(\mathrm{Index}E) E_\omega(e)$ is contructive,
%\item
%$1 - (\mathrm{Index} E)E_\omega(e) \preceq z$ in $A_\omega$,
%\item
%$A \ni x \mapsto xe \in A_\omega$ is injective.
%\end{enumerate}
%\end{dfn}
Motivated by  Theorem \ref{thm:ELPW}  we propose a trcial Rokhlin type property  in term of the tracial 2-norm  $\| \cdot\|_{2, \tau}$  for a pair $(A, \tau)$ where $A$ is a finite, infinite dimensional, simple, separable, unital $C\sp*$-algebra with a unique tracial state  $\tau$. 
\begin{defn}\label{D:probRokhlin}
Suppose that $A$ is a finite, infinite dimensional, simple, separable, unital $C\sp*$-algebra with a unique tracial state $\tau$ and Property (SP), and such that  the order on projections over $A$ is determined by the trace $\tau$. For a free ultrafilter $\omega$ let $c_{\omega, \tau} (A)=\{(x_n) \in \prod_{i=1}^{\infty}A \mid \tau(x_n^*x_n)^{1/2} \to 0 \,  \text{as}\,  n \to \omega\}$. Then set $A^{\omega}_{2, \tau} = \prod_{i=1}^{\infty} A / c_{\omega, \tau}(A)$ and call it the probabilistic ultrapower of $A$. \\
Let $P \subset A$ be an inclusion of unital $C\sp*$-algebras such that a conditional expectation $E:A \to P$ has a finite  index. We say that $E$ has the probabilistic Rokhlin property if there is a projection $e\in A^{\omega}_{2, \tau} \cap A'$ such that $(\Index E)E^{\omega}(e)=1$ where $E^{\omega}: A^{\omega}_{2, \tau} \to P^{\omega}_{2, \tau}$ is the induced map.
 \end{defn}
The following lemma is natural which is a local approximation version for Definition \ref{D:probRokhlin}.  
\begin{lem}\label{L:probRokhlin} 
Let $A$ be a simple, separable, unital, tracial  rank zero $C^*$-algebra with a unique tracial state $\tau$ and 
$P \subset A$ be an inclusion of unital  C*-algebras of index-finite type and $E: B \to A$ be the associated conditional expectation. Then the following are equivalent.
\begin{enumerate}
\item
$E$ has the probabilistic Rokhlin property .
\item
For any finite set $F \subset A$ and every $\epsilon > 0$ there exists a projection $e \in A$ such that 
\begin{itemize}
\item[(i)]
$\|ea - ae\|_{2,\tau} < \epsilon$ for any $a \in F$,
\item[(ii)]
$\|1 - (\mathrm{Index}E)E(e)\|_{2,\tau} < \epsilon$.
\end{itemize}
\end{enumerate}
\end{lem}

\begin{rem}
Under the same assumption in Lemma \ref{L:probRokhlin},  in view of Theorem \ref{thm:ELPW} we expect that $E$ has the tracial Rokhlin property if and only if $E$ has the probabilistic Rokhlin property.  But we do not know whether this holds. 
\end{rem}

When a unital $C\sp*$-algebra $A$ is equipped with  a unique tracial state $\tau$, then we denote  by $\pi_{\tau}$ the GNS-representation $A$ on $H_{\tau}$. Although not essential to the main theme,  the following result is worth noting in its own right.

\begin{prop}
Let $G$ be a finite group and $A$ a simple separable unital tracial topological rank zero $C^*$-algebra with a unique tracial state $\tau$, $\alpha$ an action of $G$ on $A$ and $E_{\alpha} \colon A \rightarrow A^\alpha$ a canonical conditional expectation. Then $\alpha$ has the tracial Rokhlin property if and only if $\overline{E_{\alpha}}^{weak}: \overline{\pi_\tau(A)}^{weak} (=M) \rightarrow \overline{\pi_\tau(A^G)}^{weak} (=N)$ has the  probabilistic Rokhlin property  in the sense that there exists an projection $e \in M^\omega \cap M'$ such that $(\overline{E_{\alpha}}^{weak})^{\omega}(e) = \frac{1}{|G|}I$.
\end{prop}

\vskip 2mm

\begin{proof}
Let $\{a_i; i \in \N\}$ be a dense subset of $A$.

Suppose that $\alpha$ has the tracial Rokhlin property. 
Then, for any $n \in \N$ and a finite set $S_n = \{a_1, a_2, \dots, a_n\}$ there exist orthogonal projection $\{e_{g,n}\}$ of $A$ such that
\begin{enumerate}
\item
$\|\alpha_g(e_{h,n}) - e_{gh,n}\|_{2,\tau} < \frac{1}{n}$ for all $g, h \in G$.
\item
$\|[e_{g,n}, a]\|_{2,\tau} < \frac{1}{n}$ for all $g \in G$, all $a \in S_n$.
\item
$\sum_{g \in G}e_{g,n} = I$.
\end{enumerate}

For any $g \in G$ set $e_g =[(\pi_\tau(e_{g,n}))]$. Then, we have 
\begin{enumerate}
\item
$\|\tilde{\alpha}_g(e_h) -e_{gh}\|_2 = \lim_{n \rightarrow \omega}\|\alpha_g(e_{h,n}) - e_{gh,n}\|_{2,\tau} = 0$ for all $g, h \in G$, that is, 
 with $\tilde{\alpha}_g(e_h) = e_{gh}$.for all $g, h \in G$.
\item
$\|[e_{g}, \pi_\tau(a)]\|_2 = \lim_{n \rightarrow \omega}\|[e_{g,n}, a]\|_{2,\tau} = 0$ for all $g \in G$, all $a \in A$, 
that is, $e_g \in M^\omega \cap M'$ for all $g \in G$.
\item
$\sum_{g \in G}e_{g} = I$,
\end{enumerate}
where $\tilde{\alpha}:G \rightarrow \mathrm{Aut}(M)$ is a canonical extended  automorphism of $\alpha$ on $A$. 

Hence, if we set $e = e_{1_G}$, then $e \in M^\omega \cap M'$ and 
$$
(\overline{E_{\alpha}}^{weak})^{\omega}(e) = \frac{1}{|G|}\sum_{g\in G}\tilde{\alpha}_g(e_{1_G}) = \frac{1}{|G|}\sum_{g\in G}e_g = \frac{I}{|G|}. 
$$

\bigskip

Conversely, suppose that $\overline{E_{\alpha}}^{weak}$ has the probabilistic Rokhlin property. Then, there exists a projection $e \in M^\omega \cap M'$ such that $\sum_{g\in G}\tilde{\alpha}_g(e) = I$. Note that $\tilde{\alpha}_g(e) \in M^\omega \cap M'$  for all $g \in G$.

Let $S \subset A$ be a finite set and $\epsilon>0$ (we may assume that any element $\|a\| \leq 1$ for all $a \in S$). We will do the same argument in the proof of \cite[Theorem 5.3]{ELPW 2010} to construct of orthogonal projections $\{e_g\}_{g \in G}$ which satisfies conditions in Theorem 5.4.

First, we obtain  $\delta_1 >0$  by applying \cite[Lemma 2.9]{Osaka-Phillips 2006} to  $\dfrac{\epsilon}{16|G|}$ in place of $\epsilon$ and $n=|G|$. Set $T = \bigcup_{g \in G}\alpha_g(S)$ and $\displaystyle \epsilon_1 = \min(\frac{1}{9|G|} \epsilon, \frac{1}{24} \epsilon, \frac{1}{9} \frac{\epsilon^2}{\delta_1}, \frac{1}{8})$. Apply \cite[Lemma 2.13]{Osaka-Phillips 2006} with $\epsilon_1$ in place of $\epsilon$, with $T$ in place of $F$, and with $G$ in place of $S$. Then we obtain projections $q, q_0 \in A$, unital finite dimensional subalgebras $E \subset qAq$ and $E_0 \subset q_0Aq_0$, and automorphisms $\phi_g \in \Aut(A)$ for $g \in G$, such that
\begin{enumerate}
    \item $\phi_1 = id_A$ and $\|\phi_g - \alpha_g\| < \epsilon_1$ for all $g \in G$,
    \item For every $g \in G$ and $x \in E$, we have $q_0\phi_g(x) = \phi_g(x)q_0$ and $q_0\phi_g(x)q_0 \in E_0$.
    \item For every $a \in T$, we have $\|qa - aq\| < \epsilon_1$ and $dist(qaq,E) < \epsilon_1$.
    \item $\tau(1 - q), \tau(1 - q_0) < \epsilon_1$.
\end{enumerate}
Apply \cite[Lemma 2.12]{Osaka-Phillips 2006} with $\epsilon_1$ in place of $\epsilon$ and $E_0 + \mathbb{C}(1-q_0)$ in place of $E$, obtaining $\delta>0$, which satisfies $\delta \leq \epsilon_1$.

Since $e \in M^\omega \cap M'$, applying \cite[Lemma 2.15]{Osaka-Phillips 2006}   and  \cite[Lemma 2.12]{Osaka-Phillips 2006} in order there exists a projection $e_0 \in A$ such that $\|\alpha_g(e_0)x - x\alpha_g(e_0)\|_{2,\tau} < \delta$ for all $x \in E_0 + \mathbb{C}(1-q_0)$ and for all  $g \in G$, and $\|\alpha_g(e_0)\alpha_h(e_0)\|_{2,\tau} < \epsilon_1$ for $g,h \in G$ with $g \neq h$. 

Let $B_0 = A \cap [E_0 + \mathbb{C}(1-q_0)]'$, the subalgebra of $A$ consisting of all elements which commute with everything in $E_0 + \mathbb{C}(1-q_0)$. Apply the choice of $\delta$ using \cite[Lemma 2.12]{Osaka-Phillips 2006} to $e_0$, obtaining a projection $f \in B_0$ which satisfies $\|f - e_0\|_{2,\tau} < \epsilon_1$. Since $q_0$ is the center of $E_0 + \mathbb{C}(1-q_0)$, the element $f_1 = fq_0$ is also a projection in $B_0$. For $g \in G$, since $q_0\phi_{g^{-1}}(E)q_0 \subset E_0$, $f_1$ commute with all elements of $q_0\phi_{g^{-1}}(E)q_0$ and with all elements of $\phi_{g^{-1}}(E)$. Therefore, $f_g = \phi_g(f_1)$ commute with all elements of $E$, including $q$, that is, $f_g$ also commutes with $1-q$.

We have the following estimation from conditions (1)-(4).
\begin{itemize}
    \item[(i)] For any $g \in G$
    \begin{align*}
        \|f_g - \alpha_g(e_0)\|_{2,\tau} &= \|\phi_g(f_1) - \alpha_g(e_0)\|_{2,\tau} \\
        &\leq \|\phi_g(f_1) - \alpha_g(f_1) + \alpha_g(f_1) - \alpha_g(e_0)\|_{2,\tau} \\
        &\leq \|\phi_g(f_1) - \alpha_g(f_1)\|_{2,\tau} + \|f_1 - e_0\|_{2,\tau} \\
        &< \epsilon_1 + \|fq_0 - e_0\|_{2,\tau} \\
        &< \epsilon_1 + \|f - e_0 - f(1-q_0)\|_{2,\tau} \\
        &\leq \epsilon_1 + \|f - e_0\|_{2,\tau} + \|f(1-q_0)\|_{2,\tau} < 3\epsilon_1.
    \end{align*}
    \item[(ii)] For any $g,h \in G$ with $g \neq h$
    \begin{align*}
        \|f_g f_h\|_{2,\tau} &\leq \|f_g f_h - \alpha_g(e_0)\alpha_h(e_0)\|_{2,\tau} + \|\alpha_g(e_0)\alpha_h(e_0)\|_{2,\tau} \\
        &= \|(f_g - \alpha_g(e_0))f_h + \alpha_g(e_0)f_h - \alpha_g(e_0)\alpha_h(e_0)\|_{2,\tau} + \epsilon_1 \\
        &\leq \|f_g - \alpha_g(e_0)\|_{2,\tau} + \|f_h - \alpha_h(e_0)\|_{2,\tau} + \epsilon_1 < 7\epsilon_1 < \delta_1.
    \end{align*}
    \item[(iii)] Since $\|\phi_g(f_1) - \alpha_g(f_1)\|_{2,\tau} < \epsilon_1$, for any $g,h \in G$
    \begin{align*}
        \|f_{gh} - \alpha_g(f_h)\|_{2,\tau} &\leq \|f_{gh} - \alpha_{gh}(f_1)\|_{2,\tau} + \|\alpha_{gh}(f_1) - \alpha_g(f_h)\|_{2,\tau} \\
        &< \epsilon_1 + \|\alpha_g(\alpha_h(f_1) - f_h)\|_{2,\tau} < \epsilon_1 + \epsilon_1 = 2\epsilon_1.
    \end{align*}
    \item[(iv)] Since for any $a \in T$ dist($qaq,E) < \epsilon_1$, that is, there exists $x_a \in E$ such that $\|qaq - x_a\|_{2,\tau} < \epsilon_1$, for any $g \in G$ and $a \in T$
    \begin{align*}
        \|f_g a - a f_g\|_{2,\tau} &= \|f_g qaq - qaqf_g + f_gqa(1-q) + f_g(1-q)aq + f_g(1-q)a(1-q) \\
        &\quad + qa(1-q)f_g + (1-q)aqf_g + (1-q)a(1-q)f_g\|_{2,\tau} \\
        &\leq \|f_g(qaq - x_a) + f_g x_a - x_a f_g - (qaq - x_a)f_g\|_{2,\tau} + 6\epsilon_1 \\
        &< 2\epsilon_1 + 6\epsilon_1 = 8\epsilon_1 \quad (\text{note that } f_g x_a = x_a f_g).
    \end{align*}
    \item[(v)]
    $\tau(1 - \sum_{g \in G} f_g) = \tau(1 - \sum_{g \in G} \alpha_g(e_0)) + \tau(\sum_{g \in G} (\alpha_g(e_0) - f_g)) < \epsilon_1 + |G|\epsilon_1$.
\end{itemize}
From the condition (ii), we have orthogonal projections $\{p_g\}_{g \in G}$ such that $\|p_g - f_g\|_{2,\tau} < \frac{\epsilon}{16}$. Hence, for any $g,h \in G$ with $g \neq h$
\begin{align*}
    \|p_{gh} - \alpha_g(p_h)\|_{2,\tau} &= \|p_{gh} - f_{gh} + f_{gh} - \alpha_g(f_h) + \alpha_g(f_h) - \alpha_g(p_h)\|_{2,\tau} \\
    &= \|p_{gh} - f_{gh}\|_{2,\tau} + \|f_{gh} - \alpha_g(f_h)\|_{2,\tau} + \|\alpha_g(f_h) - \alpha_g(p_h)\|_{2,\tau} \\
    &< \frac{\epsilon}{16} + \|f_{gh} - \alpha_g(f_h)\|_{2,\tau} + \|f_h - p_h\|_{2,\tau} < \frac{\epsilon}{16} + 2\epsilon_1 + \frac{\epsilon}{16} < \epsilon.
\end{align*}
For any $g \in G$ and $a \in T$
\begin{align*}
    \|p_g a - a p_g\|_{2,\tau} &= \|(p_g - f_g)a + f_g a - a f_g - a(p_g - f_g)\|_{2,\tau} \\
    &< \|p_g - f_g\|_{2,\tau} + \|f_g a - a f_g\|_{2,\tau} + \|p_g - f_g\|_{2,\tau} \\
    &< \frac{\epsilon}{16} + 8\epsilon_1 + \frac{\epsilon}{16} < \frac{\epsilon}{8} + \frac{\epsilon}{24} < \epsilon.
\end{align*}
$\tau(1 - \sum_{g \in G} p_g) = \tau(1 - \sum_{g \in G} f_g) + \sum_{g \in G} \tau(f_g - p_g) < \epsilon_1 + |G|\epsilon_1 + \frac{\epsilon}{16|G|}|G| < \frac{\epsilon}{24} + \frac{\epsilon}{9} + \frac{\epsilon}{16} < \epsilon$.

Therefore, for any finite set $S$ and $\epsilon>0$ there exist orthogonal projection $p_g \in A$ for $g \in G$ satisfying conditions (1),(2) and (3') in \cite[Theorem 5.3]{ELPW 2010}. Hence, this implies that for any finite set $S \subset A$ and $\epsilon>0$ there exists orthogonal projections $\{e_g\}_{g \in G}$ satisfying conditions (1),(2), and (3) in Theorem 5.4.
Hence, we conclude that an action $\alpha$ has the tracial Rokhlin property.
\end{proof}

%%%%%%%%%%%%%%%%%%%%%%%%%%%%%%%%%%%%%%%%%%%%%%%%%%%%%%%%%%%%%%%%%%%%%%%%%%%%%%%%%%%%%%%%%%%%
\vskip 2mm
\begin{prop}\label{prp:Tracial Rokhlin}
Let $A \subset B$ be an inclusion of simple, unital tracial topological rank zero, C*-algebras with a unique tracial state $\tau$ and a conditional expectation $E: B \rightarrow A$ is of index-finite type. Let $M$ ( resp. $N$) be the weak closure of the GNS construction $\pi_\tau(B)$ of $B$ (resp. the weak closure of the GNS construction $\pi_\tau(A)$). Then, if there is a projection $e \in M^\omega \cap M'$ such that $E^\omega(e) = \frac{1}{[M:N]}$, where $E^\omega\colon M^\omega \rightarrow N^\omega$ be of index-finite type, $E$ has the  probabilistic Rokhlin property.
\end{prop}
\begin{proof}
We may identify $A$ (resp. $B$) with $\pi_\tau(A)$ (resp. $\pi_\tau(B)$). It is enough to show that (1), (2) in Lemma \ref{L:probRokhlin}; Let $F \subset \pi_\tau(B)$ a finite set  and $\epsilon>0$. Since $\pi_\tau(B)$ has tracial rank zero, there is a projection $q \in  \pi_\tau(B)$ and a unital finite dimensional $D \subset q\pi_\tau(B)q$ such that $\|qa -aq\| < \epsilon$ $(a \in F)$, $\mathrm{dist}(qaq, D) < \epsilon$, and $\tau(1 - q) < \epsilon$. 

By the assumption,  there is $e=[(e_n)_n] \in M_\omega \cap M'$ such that $[M:N] (\overline{E}^{weak})^{\omega}(e)=1$.  By applying \cite[Lemma 2.12]{Osaka-Phillips 2006} to $D$, $\epsilon>0$ we obtain $\delta >0$. Therefore,  there exists $e_n \in M$ such that  $\|E(e_n) - \frac{1}{[M:N]}\|_{2, \tau} < \frac{\epsilon}{[M:N]}$ and  $\|e_n x - xe_n\|_{2,\tau} < \delta$  for all $x$ in the matrix system in $D$.   May assume that $e_n \in \pi_\tau(B)$ by \cite[Lemma 2.15]{Osaka-Phillips 2006}. 
Then, there is a projection $p \in \textcolor{red}{\pi_\tau(B)}$ such that $p$ commutes with $D$ and $\|e_n - p\|_{2,\tau} < \epsilon$.

Note that $\mathrm{Index}E = [M:N]$ by \cite[Corollary~3.7]{MI 2002}. It follows that \[\|1 - \mathrm{Index}EE(e_n)\|_{2,\tau} < \epsilon\].
 
Since $pq = qp$, $qp$ is projection. %We show that For any $a \in F$ \textcolor{red}{$[a, qp]$} is small. 
Take $d \in D$ such that $\|qaq - d\| < \epsilon$. Then,

%\begin{align*}
%\|dpq - pqd\| &= \|(d-qaq)pq + qaqpq -pq(d-qaq) -pqqaq\|\\
%&\leq \|d-qaq\| + \|qapq - pqaq\| + \|d - qaq\|\\
%&< 2\epsilon + \|qapq - p
%\end{align*}

\begin{align*}
\|a - (d + (1-q)a(1-q))\| &= \|qaq - d + [q,a](1-q) + (1-q)[q,a]\|\\
&\leq \|qaq -d\| + 2\|[q,a]\| < 3\epsilon
\end{align*}
Hence,
\begin{align*}
\|(qp)a - a(qp)\| &= \|qp\{a-(d + (1-q)a(1-q))\} + qpd - \{a-(d + (1-q)a(1-q))\})qp -dqp\|\\
&\leq 2\|a-(d + (1-q)a(1-q))\|  + \|[d, qp]\| < 6\epsilon \ (qd = d, pd = dp)\  
\end{align*} for all $a\in F$.

%Since $|\tau(p-e_n)|^2 \leq \tau((p-e_n)^2) = \|p-e_n\|_\tau^2 < \epsilon^2$, $|\tau(p-e_n)| < \epsilon$. Hence, 

Since $\|p -e_n\|_{2, \tau} < \epsilon$, 

\begin{align*}
\|pq - e_n\|_{2,\tau} &= \|p(q-1) + p -e_n\|_{2,\tau}\\
&\leq \|p(q-1)\|_{2,\tau} + \|p -e_n\|_{2,\tau}\\
&\leq \tau(1-q)^{1/2} + \epsilon\\
&< \sqrt{\epsilon} + \epsilon.
\end{align*}

%\begin{align*}
%|\tau (pq - e_n)| &= |\tau(p(q-1) + p-e_n)|\\
%&\leq |\tau(p(q-1))| + |\tau(p-e_n)|\\
%&\leq \tau(p)^\frac{1}{2}\tau(1-q)^\frac{1}{2} +\epsilon\\
%&< \sqrt{\epsilon} + \epsilon.
%\end{align*}

%\begin{align*}
%|\tau(pq -e)| &= |\tau(p(q-1) + p-e)|\\
%&\leq |\tau(q-e)| + |\tau((p-e)e)|\\
%&\leq |\tau(q-e)| + |\tau(p-e) |< 2\epsilon.
%\end{align*}

Therefore, for taking sufficiently large $n$

\begin{align*}
\|1 - \mathrm{Index}E E(pq)\|_{2,\tau}
&\leq \|1 - \mathrm{Index}E E(e_n)\|_{2,\tau} + \mathrm{Index}E\|(E(e_n) - E(pq))\|_{2,\tau}\\
&<\epsilon + \mathrm{Index}E\|E(e_n-pq)\|_{2,\tau}\\
&\leq \epsilon + \mathrm{Index}E\|e_n - pq\|_{2,\tau} \ (\tau \circ E = \tau\hbox{})\\
&< \epsilon + \mathrm{Index}E (\epsilon + \sqrt{\epsilon})
\end{align*}

%\begin{align*}
%|\tau(1 - \mathrm{Index}E E(pq))|
%&\leq |\tau(1 - \mathrm{Index}E E(e_n)| + |\mathrm{Index}E\textcolor{red}{\tau}(E(e_n) - E(pq))|\\
%&< \epsilon + \mathrm{Index}E\textcolor{red}{ |\tau(e_n-pq)|} \\
%&\textcolor{red}{< \epsilon + \mathrm{Index}E (\epsilon + \sqrt{\epsilon})}\\
%= (1 + \mathrm{Index}(E))\epsilon.
%\end{align*}

%\color{red}
%We have to reconsider the below argument.

%\color{black}

%From Lemma~\ref{lem:perturbation} and \textcolor{red}{Kaplansky's density theorem} we can replace $pq$ by a contractive positive element $r \in M$ such that $\|ra - ar\|$ $\forall a \in F)$, $\mathrm{Index}E E(r)$ is positive contructive, and $\tau(1 - \mathrm{Index}E E(r))$ are small. 
Therefore, by Lemma~\ref{L:probRokhlin},  $E$ has the probabilistic  Rohlin property.
\end{proof}

Our final piece in the construction of irreducible inclusion of finite $C\sp*$-algebras with the probabilistic  Rokhlin property come from the following famous result.   
\begin{prop}\label{prp:PP}\cite[Theorem~4.1.2]{Popa 1994}
Let $N \subset M$ of  strongly amenable, injective $\rm II_1$-factor of index finite type.
there exist simple AF C*-algebras $A \subset B$ of index-finite type such that 
$\overline{A}^{weak} = N$ and $\overline{B}^{weak} = M$.
\end{prop}
\vskip 2mm

Now we are finally ready to present our main theorem of this section. 

\begin{thm}\label{thm:tracial rikhlin property}
There is an inclusion $A \subset B$ of simple AF-C*-algebras which has the probabilistic Rokhlin property with non-integer index.
\end{thm}

\begin{proof}
From Proposition~\ref{prp:PP},  if $N \subset M$ is  an  inclusion of  strongly amenable, injective $\rm II_1$-factors of index-finite type, then 
there exist simple AF $C^*$-algebras $A \subset B$ of index-finite type such that 
$\overline{A}^{weak} = N$ and $\overline{B}^{weak} = M$.
Hence, from Conjecture~\ref{cnj:rokhlin} and Proposition~\ref{prp:Tracial Rokhlin} we know that the inclusion $A \subset B$ has the  tracial Rokhlin property. Moreover, its index can be non-integer.
\end{proof}

In summary, we have provided constructions of previously unknown examples of irreducible Rokhlin inclusions with non-integer indices, exhibiting the non-classical quantum symmetries arising in subfactor theory and encoded within the structure of unitary tensor categories. A natural and important direction for future research is to investigate more systematically whether integer-index inclusions can genuinely depart from classical fixed-pint algebra frameworks. Such efforts promise deeper insights into the subtle boundary between classical and quantum symmetries in operator algebras. Related to this, we list important questions in our opinion.   

\begin{question}
Is there an inclusion $A \subset B$ of finite simple unital C*-algebras with the Rokhlin property such that its index is non-integer?
\end{question}

A unital simple purely infinite $C^*$-algebra $A$ is said to be in the Cuntz standard form if $[1] = 0$ holds in $K_0(A)$. That is, there is an inclusion from $\mathcal{O}_2$ to $A$. 
\begin{question}
Is there an irreducible inclusion $A \subset B$ of unital simple purely infinite $C^*$-algebras (which does not isomorphic to $\mathcal{O}_2$) in the Cuntz standard form with a conditional expectation $E$ from $B$ onto $A$, whose index is  a non-integer, such that it has the Rokhlin property? 
\end{question}

\begin{question}
If the inclusion $P \subset A$ of index-finite type has the Rokhlin property and the associated conditional expectation $E$ has an integer index, say $n$, then what kind of conditions do we need to see that there are  a finite group $G$ with $|G|=n$ and an action $\alpha: G \curvearrowright A$ such that $A^{\alpha}=P$ and $E=E_{\alpha}$. 
\end{question}
\vskip 2mm

\section{Acknowledgement}
A part of the first author's research was done during his visit to Korea Institute of Advanced Study in the year of 2024. He would like to aknowledge the support from the institute.   

\end{document}